\documentclass[11pt]{article}
\usepackage[T1]{fontenc}
\usepackage{lmodern}
\usepackage{amsmath,amssymb,amsthm}
\usepackage{microtype}
\usepackage[margin=1.07in]{geometry}
\usepackage{titling}
\usepackage{hyperref}
\hypersetup{colorlinks=true,linkcolor=black,citecolor=black,urlcolor=blue}

\newtheorem{theorem}{Theorem}
\newtheorem{lemma}[theorem]{Lemma}
\newtheorem{corollary}[theorem]{Corollary}
\theoremstyle{remark}
\newtheorem{remark}[theorem]{Remark}

\newcommand{\Ric}{\operatorname{Ric}}
\newcommand{\C}{\mathbb C}

\title{K\"ahler metrics with varying cone angles and Ricci curvature}
\author{Martin de Borbon\\[4pt]
\small Department of Mathematical Sciences\\
\small Loughborough University\\
\small Loughborough, LE11 3TU, United Kingdom\\
\small \texttt{m.de-borbon@lboro.ac.uk}}
\date{}

\begin{document}
\maketitle

\begin{abstract}
We prove that the cone angle of a K\"ahler metric is constant if its Ricci curvature is locally in \(L^2\). In particular, a K\"ahler--Einstein metric cannot have a varying cone angle. The proof is elementary and assumes only uniform comparison with the corresponding model.
\end{abstract}



\section{Introduction}

Atiyah and LeBrun asked in \cite{AtiyahLeBrun} whether an Einstein metric can have a cone angle varying along a codimension-two submanifold. We give a negative answer in the K\"ahler setting. More generally, we show that the cone angle is constant if the Ricci curvature is locally square-integrable.

We work on \(\C^n\), \(n\geq2\), with coordinates \((z,w)\), where \(z\in\C\)  and \(w=(w_1,\ldots,w_{n-1}) \in \C^{n-1}\).
Let \(\Delta^* = \{0<|z|<1\} \subset \C\) and let \(B = \{|w| < 1\} \subset \C^{n-1}\).
For a smooth function \(\beta:B\to(0,\infty)\), set
\begin{equation}\label{eq:model}
	g_\beta=|z|^{2\beta(w)-2}|dz|^2+\sum_{j=1}^{n-1}|dw_j|^2.
\end{equation}
This is a Hermitian metric on \(\Delta^*\times B\) whose restriction to each punctured disc \(\Delta^*\times\{w\}\) is the standard flat metric with cone angle \(2\pi\beta(w)\) at \(z=0\).

\begin{theorem}\label{thm:main}
	Let \(g\) be a K\"ahler metric on \(\Delta^*\times B\) with 
	\begin{equation}\label{eq:comparison}
		C_0^{-1}g_\beta\leq g\leq C_0g_\beta
	\end{equation}
	for some constant \(C_0\geq1\).
	If
	\begin{equation}\label{eq:L2-ricci-bound}
		\int_{\Delta^* \times B} |\Ric(g)|_g^2 \, dV_g < \infty ,
	\end{equation}
	then \(\beta\) is constant on \(B\).
\end{theorem}

\begin{corollary}\label{cor:bounded}
	Let \(X\) be a complex manifold and \(D\subset X\) a connected smooth divisor. Let \(g\) be a K\"ahler metric on \(X\setminus D\) and \(\beta:D\to(0,\infty)\) a smooth function. Suppose that near every point of \(D\) there are holomorphic coordinates \((z,w)\) in which \(D=\{z=0\}\) and the comparison \eqref{eq:comparison} holds with \(\beta(w)=\beta(0,w)\). If \(|\Ric(g)|_g\) is bounded near \(D\), then \(\beta\) is constant. In particular, a K\"ahler--Einstein metric with these conical asymptotics has constant cone angle.
\end{corollary}

\begin{proof}
	The model \eqref{eq:model} has finite volume near \(D\). Thus bounded Ricci curvature implies local \(L^2\) integrability. The corollary follows from Theorem~\ref{thm:main} and the connectedness of \(D\).
\end{proof}

\begin{remark}
	The comparison condition \eqref{eq:comparison} is independent of the choice of holomorphic coordinates adapted to the divisor when \(\beta\leq1\). This is no longer true for \(\beta>1\), as can be seen by pulling back the model metric with constant \(\beta>1\) under the map \((z,w)\mapsto(z,z+w)\).
\end{remark}

\begin{remark}\label{rem:sharpness}
	The exponent two in Theorem~\ref{thm:main} is optimal. The potential
	\(\varphi(z,w)=|w|^2+|z|^{2\beta(w)}\)
	defines, locally near \(\{z=0\}\), a K\"ahler metric satisfying \eqref{eq:comparison}. A direct calculation shows that the Ricci curvature is locally in \(L^p\) for every \(p<2\). Moreover, the scalar curvature is bounded if \(\beta^{-1}\) is harmonic. Thus bounded scalar curvature does not imply constant cone angle. 
	
	Oliveira and Sena-Dias \cite{OliveiraSenaDias} constructed scalar-flat K\"ahler metrics on complex surfaces with varying cone angles along a smooth complex curve.
\end{remark}

\begin{remark}\label{rem:riemannian}
	Let \((\Sigma,h)\) be a Riemannian manifold and \(\beta:\Sigma\to(0,\infty)\) a smooth function. For the metric \(g_0=d\rho^2+\beta(x)^2\rho^2\,d\theta^2+h\),
	we have \(\Ric(g_0)(\partial_\rho,V)=-\rho^{-1}d\log\beta(V)\) for every vector \(V\) tangent to \(\Sigma\). Near a point where \(d\beta\neq0\), the norm of the Ricci curvature blows up like \(\rho^{-1}\).
\end{remark}

\section{Proof of the theorem}

The idea of the proof is to recover the angle function from the logarithm of the determinant of the metric. Since the comparison condition gives no control over derivatives, we take circular averages in the normal variable. The mixed components of the Ricci tensor then control the variation of these averages along the divisor.

Write \(z=se^{i\theta}\), where \(s=|z|\), and set \(t=-\log s\). Let \(u=\log\det(g_{i\bar j})\) and consider 
\begin{equation}\label{eq:averages}
	m(t,w)=\frac1{2\pi}\int_0^{2\pi}u(se^{i\theta},w)\,d\theta.
\end{equation}

\begin{lemma}\label{lem:average}
Let \(a=\frac12 \partial_t m\). Then
\begin{equation}\label{eq:contour}
 \partial_{\bar w_j}a(t,w)
 =\frac{1}{2\pi i}\int_{|z|=s}\Ric_{z\bar w_j}\,dz.
\end{equation}
\end{lemma}

\begin{proof}
On the circle \(|z|=s\), we have \(\partial_z u\,dz=\frac12(\partial_\theta u-i\partial_t u)\,d\theta\). Since \(\int_0^{2\pi}\partial_\theta u\,d\theta=0\), we obtain
\[
 a(t,w)=-\frac{1}{2\pi i}\int_{|z|=s}\partial_z u\,dz.
\]
Differentiating under the integral sign and using \(\Ric_{z\bar w_j}=-\partial_z\partial_{\bar w_j}u\) gives \eqref{eq:contour}.
\end{proof}

Throughout the rest of the proof, \(C\) denotes a positive constant, depending only on \(C_0\) and \(n\), which may change from line to line. For positive functions or volume densities, we write \(f_1\asymp f_2\) if \(C^{-1}f_2\leq f_1\leq Cf_2\) for some constant \(C\geq1\).

\begin{lemma}\label{lem:energy}
Under the comparison \eqref{eq:comparison},
\begin{equation}\label{eq:energy}
 \int_{0}^{\infty}\!\int_B|\nabla_w a(t,w)|^2\,dw\,dt
 \leq C\int_{\Delta^*\times B}|\Ric(g)|_g^2\,dV_g,
\end{equation}
where \(\nabla_w\) and \(dw\) are the Euclidean gradient and volume measure in the \(w\)-variables.
\end{lemma}

\begin{proof}
Writing \(z=se^{i\theta}\) in \eqref{eq:contour}, so that \(dz=ise^{i\theta}\,d\theta\), and applying Cauchy--Schwarz, we obtain
\[
 |\partial_{\bar w_j}a(t,w)|^2
 \leq\frac{s^2}{2\pi}\int_0^{2\pi}
 |\Ric_{z\bar w_j}(se^{i\theta},w)|^2\,d\theta.
\]
Since \(a\) is real-valued, we have \(|\nabla_w a|^2=4\sum_{j=1}^{n-1}|\partial_{\bar w_j}a|^2\). Summing the above inequalities over \(j\),
\begin{equation}\label{eq:circlebound}
 |\nabla_w a(t,w)|^2
 \leq\frac{2s^2}{\pi}\int_0^{2\pi}\sum_{j=1}^{n-1}
 |\Ric_{z\bar w_j}(se^{i\theta},w)|^2\,d\theta.
\end{equation}
The comparison \eqref{eq:comparison} implies that the vectors \(s^{1-\beta(w)}\partial_z\) and \(\partial_{w_j}\) have uniformly bounded \(g\)-lengths. Hence
\begin{equation}\label{eq:riccomponent}
 \sum_{j=1}^{n-1}|\Ric_{z\bar w_j}|^2
 \leq Cs^{2\beta(w)-2}|\Ric(g)|_g^2.
\end{equation}
Combining \eqref{eq:circlebound} and \eqref{eq:riccomponent}, the radial factors give \(s^2s^{2\beta(w)-2}=s^{2\beta(w)}\). On the other hand,
the metric comparison \eqref{eq:comparison} gives \(dV_g\asymp s^{2\beta(w)}\,dt\,d\theta\,dw\).  Integrating \eqref{eq:circlebound} over \((0,\infty)\times B\) proves \eqref{eq:energy}.
\end{proof}

\begin{lemma}\label{lem:growth}
Let \(m\) be a smooth function on \((0,\infty)\times B\) of the form
\begin{equation}\label{eq:growth}
 m(t,w)=2t\,b(w)+r(t,w),
\end{equation}
where \(b\) is continuous on \(B\) and \(r\) is bounded on \((0,\infty)\times B\). If \(a=\frac12\partial_t m\) satisfies
\begin{equation}\label{eq:finiteenergy}
 \int_0^\infty\!\int_B|\nabla_w a(t,w)|^2\,dw\,dt<\infty,
\end{equation}
then \(b\) is constant on \(B\).
\end{lemma}

\begin{proof}
Fix \(t_0>0\) and set
\[
 h_t(w)=\frac{m(t,w)-m(t_0,w)}{2(t-t_0)},\qquad t>t_0.
\]
By \eqref{eq:growth},
\[
 h_t(w)-b(w)=\frac{r(t,w)-r(t_0,w)}{2(t-t_0)}.
\]
The boundedness of \(r\) implies that \(h_t\to b\) uniformly on \(B\) as \(t\to\infty\). Since \(\partial_t m=2a\), we have
\[
 \nabla_w h_t(w)=\frac{1}{t-t_0}\int_{t_0}^t\nabla_w a(\tau,w)\,d\tau.
\]
For each fixed \(w\), Cauchy--Schwarz gives
\[
 \left|\int_{t_0}^t\nabla_w a(\tau,w)\,d\tau\right|^2
 \leq(t-t_0)\int_{t_0}^t|\nabla_w a(\tau,w)|^2\,d\tau.
\]
Integrating over \(B\), we obtain
\[
 \|\nabla_w h_t\|_{L^2(B)}^2
 \leq\frac{1}{t-t_0}
 \int_{t_0}^\infty\!\int_B|\nabla_w a(\tau,w)|^2\,dw\,d\tau
 \longrightarrow0.
\]
Let \(B'\Subset B\) be any ball and denote the average of a function \(f\) on \(B'\) by \(\overline f=|B'|^{-1}\int_{B'}f(w)\,dw\). The Poincar\'e inequality gives
\[
 \|h_t-\overline{h_t}\|_{L^2(B')}
 \leq C_{B'}\|\nabla_w h_t\|_{L^2(B')}
 \longrightarrow 0.
\]
On the other hand, since \(h_t\to b\) uniformly on \(B'\), we have \(\overline{h_t}\to\overline b\). It follows that \(h_t-\overline{h_t}\to b-\overline b\) in \(L^2(B')\). Hence
\[
 \|b-\overline b\|_{L^2(B')}=0.
\]
Since \(B'\) was arbitrary, \(b(w)\) is constant on \(B\).
\end{proof}

\begin{proof}[Proof of Theorem~\ref{thm:main}]
Taking determinants in \eqref{eq:comparison}, we obtain \(\det(g_{i\bar j})\asymp s^{2\beta(w)-2}\). Consequently,
\begin{equation}\label{eq:det}
 u=2t(1-\beta(w))+q(z,w),
\end{equation}
where \(|q|\leq C\) on \(\Delta^*\times B\). Taking the circular average gives \eqref{eq:growth} with \(b(w)=1-\beta(w)\) and \(r\) the circular average of \(q\), which is bounded. By Lemma~\ref{lem:energy}, the finite Ricci energy \eqref{eq:L2-ricci-bound} implies \eqref{eq:finiteenergy}. Thus Lemma~\ref{lem:growth} shows that \(\beta(w)\) is constant on \(B\).
\end{proof}

\section*{Statements and Declarations}

\noindent\textbf{Funding.}
No specific funding was received for this work.

\medskip
\noindent\textbf{Competing interests.}
The author declares no competing interests.

\medskip
\noindent\textbf{Data availability.}
No datasets were generated or analysed in this study.

\medskip
\noindent\textbf{Use of AI.}
ChatGPT was used to assist with mathematical discussions,
drafting and editing of the manuscript. The author takes
full responsibility for the final version.

\end{document}